\documentclass[12pt]{article}
\usepackage{graphicx}
\usepackage{amsmath,amsthm,amssymb,enumerate}
\usepackage{euscript,mathrsfs}
\usepackage{color}
\usepackage{dsfont}
\usepackage[left=2cm,right=2cm,top=3.5cm,bottom=3.5cm]{geometry}
\usepackage{color}
\usepackage[framemethod=tikz]{mdframed}
\allowdisplaybreaks

\usepackage{soul}
\usepackage{esint}
\catcode`\@=11 \@addtoreset{equation}{section}

\catcode`\@=12

\allowdisplaybreaks

\newtheorem{Theorem}{Theorem}[section]
\newtheorem{Proposition}[Theorem]{Proposition}
\newtheorem{Lemma}[Theorem]{Lemma}
\newtheorem{Corollary}[Theorem]{Corollary}

\theoremstyle{definition}
\newtheorem{Definition}[Theorem]{Definition}

\newtheorem{Remark}[Theorem]{Remark}

\newcommand{\bTheorem}[1]{
\begin{Theorem} \label{T#1} }
\newcommand{\eT}{\end{Theorem}}

\newcommand{\bProposition}[1]{
\begin{Proposition} \label{P#1}}
\newcommand{\eP}{\end{Proposition}}

\newcommand{\bLemma}[1]{
\begin{Lemma} \label{L#1} }
\newcommand{\eL}{\end{Lemma}}

\newcommand{\bCorollary}[1]{
\begin{Corollary} \label{C#1} }
\newcommand{\eC}{\end{Corollary}}

\newcommand{\bRemark}[1]{
\begin{Remark} \label{R#1} }
\newcommand{\eR}{\end{Remark}}

\newcommand{\bDefinition}[1]{
\begin{Definition} \label{D#1} }
\newcommand{\eD}{\end{Definition}}

\newcommand{\tvm}{\widetilde{\vc{m}}}

\newcommand{\bfphi}{\boldsymbol{\varphi}}

\newcommand{\bFormula}[1]{
\begin{equation} \label{#1}}
\newcommand{\eF}{\end{equation}}

\newcommand{\Ov}[1]{\overline{#1}}

\newcommand{\vr}{\varrho}

\newcommand{\tvr}{\tilde \vr}
\newcommand{\tvu}{{\tilde \vu}}

\newcommand{\vu}{\vc{u}}
\newcommand{\vm}{\vc{m}}

\newcommand{\vc}[1]{{\bf #1}}

\newcommand{\Div}{{\rm div}_x}
\newcommand{\Grad}{\nabla_x}

\newcommand{\dx}{\,{\rm d} {x}}

\newcommand{\dt}{\,{\rm d} t }

\newcommand{\intTd}[1]{\int_{\mathbb{T}^d} #1 \ \dx}

\newcommand{\vv}{\vc{v}}

\newcommand{\ep}{\varepsilon}

\newcommand{\br}{ \nonumber \\ }

\newcommand{\Td}{\mathbb{T}^d}

\def\softd{{\leavevmode\setbox1=\hbox{d}%
          \hbox to 1.05\wd1{d\kern-0.4ex{\char039}\hss}}}
\definecolor{Cgrey}{rgb}{0.85,0.85,0.85}
\definecolor{Cblue}{rgb}{0.50,0.85,0.85}
\definecolor{Cred}{rgb}{1,0,0}
\definecolor{fancy}{rgb}{0.10,0.85,0.10}

\newcommand\Cbox[2]{%
    \newbox\contentbox%
    \newbox\bkgdbox%
    \setbox\contentbox\hbox to \hsize{%
        \vtop{
            \kern\columnsep
            \hbox to \hsize{%
                \kern\columnsep%
                \advance\hsize by -2\columnsep%
                \setlength{\textwidth}{\hsize}%
                \vbox{
                    \parskip=\baselineskip
                    \parindent=0bp
                    #2
                }%
                \kern\columnsep%
            }%
            \kern\columnsep%
        }%
    }%
    \setbox\bkgdbox\vbox{
        \color{#1}
        \hrule width  \wd\contentbox %
               height \ht\contentbox %
               depth  \dp\contentbox
        \color{black}
    }%
    \wd\bkgdbox=0bp%
    \vbox{\hbox to \hsize{\box\bkgdbox\box\contentbox}}%
    \vskip\baselineskip%
}

\mdfdefinestyle{MyFrame}{%
	linecolor=black,
	outerlinewidth=1pt,
	roundcorner=5pt,
	innertopmargin=\baselineskip,
	innerbottommargin=\baselineskip,
	innerrightmargin=10pt,
	innerleftmargin=10pt,
	backgroundcolor=white!20!white}

\date{}

\begin{document}


\title{On the density of ``wild'' initial data for the barotropic Euler system}

\author{Elisabetta Chiodaroli 
\and Eduard Feireisl \thanks{The work of E.F. was partially supported by the
Czech Sciences Foundation (GA\v CR), Grant Agreement
21--02411S. The Institute of Mathematics of the Academy of Sciences of
the Czech Republic is supported by RVO:67985840. } 
}

\date{}

\maketitle

\medskip

\centerline{Dipartimento di Matematica, Largo Bruno Pontecorvo} 
	
\centerline{Univerit\` a di Pisa, Italy}

\medskip

\centerline{Institute of Mathematics of the Academy of Sciences of the Czech Republic}

\centerline{\v Zitn\' a 25, CZ-115 67 Praha 1, Czech Republic}

\date{}

\maketitle

\begin{abstract}
	
	We show that the set of ``wild data'', meaning the initial data for which the barotropic Euler system
	admits infinitely many \emph{admissible entropy} solutions, is dense in the $L^p-$topology of the phase space.

\end{abstract}

\bigskip


\bigskip

\section{Introduction}
\label{i}

The concept of wild data/solution appeared recently in the context of the ill--posedness results 
obtained via the method of convex integration, see e.g. Buckmaster et al \cite{BuDeSzVi}, \cite{BucVic} and the 
references cited therein. In contrast with the concept of wild solution that may be ambiguous, the wild data
can be identified with those that give rise to infinitely many solutions of a given problem on any (short) 
time interval, see Definition \ref{wD1} below. Sz\'ekelyhidi and Wiedemann \cite{SzeWie} showed that the set of wild data 
for the \emph{incompressible} Euler system in the framework of weak solutions satisfying the global energy inequality is dense in the 
$L^p-$topology of the phase space. Our goal is to discuss the problem in the class of weak entropy solutions of the barotropic Euler system 
describing the motion of a compressible fluid.

\subsection{Barotropic Euler system}

We consider the \emph{barotropic Euler system}:

\begin{align}
	\partial_t \vr + \Div (\vr \vu) &= 0, \label{i1}\\
	\partial_t (\vr \vu) + \Div (\vr \vu \otimes \vu) + \Grad p(\vr) &= 0 \label{i2}
\end{align}
describing the time evolution of the mass density $\vr = \vr(t,x)$ and the velocity $\vu = \vu(t,x)$ of a compressible inviscid fluid.
For simplicity, we impose the space periodic boundary conditions identifying the fluid domain with the flat torus: 
\begin{equation} \label{i3}
	\Td = \left( [-1,1]|_{\{ -1; 1 \} } \right)^d, \ d=2,3.
\end{equation} 
The same method can be used to obatin similar results for fluids occupying a bounded domain $\Omega \subset R^d$, with impermeable boundary 
\begin{equation} \label{i3a}
	\vu \cdot \vc{n}|_{\partial \Omega} = 0,\ \vc{n} - \mbox{the outer normal vector to}\ \partial \Omega, 
	\end{equation}
see Section \ref{cr}.

\begin{Definition}[\bf Admissible entropy solution] \label{iD1}
	
	We say that $(\vr, \vu)$ is \emph{admissible entropy solution} to the Euler system \eqref{i1}--\eqref{i3} in $(0,T) \times \Td$ with 
	initial data 
\[
\vr(0, \cdot) = \vr_0,\ \vu(0, \cdot) = \vu_0
\]
if the following holds:	

\begin{align} 
	\int_0^T \intTd{ \Big[ \vr \partial_t \varphi + \vr \vu \cdot \Grad \varphi \Big] } \dt &= - \intTd{ \vr_0 \varphi(0, \cdot) }, \label{i4} \\ 
	\mbox{for any} \ \varphi &\in C^1_c([0,T) \times \Td), \br
	\int_0^T \intTd{ \Big[ \vr \vu \cdot \partial_t \bfphi + \vr \vu \otimes \vu : \Grad \bfphi + p(\vr) \Div \bfphi \Big] } \dt &= - \intTd{ \vr_0 \vu_0 \cdot \bfphi(0, \cdot) } \label{i5}\\
	\mbox{for any} \ \bfphi &\in C^1_c([0,T) \times \Td; R^d), \nonumber
	\end{align}
with the energy inequality:
\begin{align} 
\int_0^T &\intTd{ \left[ \left( \frac{1}{2} \vr |\vu|^2 + P(\vr) \right) \partial_t \varphi 
	+ \left( \frac{1}{2} \vr |\vu|^2 + P(\vr) + p(\vr) \right) \vu \cdot \Grad \varphi   \right]} \dt 
\br
&\geq - \intTd{ \left( \frac{1}{2} \vr_0 |\vu_0|^2 + P(\vr_0) \right) \varphi (0, \cdot) } 	
\label{i6} \\ 
\mbox{for any}\ \varphi &\in C^1([0,T) \times \Td ),\ \varphi \geq 0,
\nonumber
	\end{align}
where we have introduced the
pressure potential: 
\[
P'(\vr) \vr - P(\vr) = p(\vr). 
\]

\end{Definition}

\subsection{Wild data, main result}
\label{w}

\begin{Definition}[\bf Wild data] \label{wD1}
	We say that the initial data $\vr_0, \vu_0$ are \emph{wild} if there exists $T_{\rm w} > 0$ such that 
	the Euler system admits infinitely many admissible entropy solutions $(\vr, \vu)$ on any interval $[0,T]$, $0 < T < T_{\rm w}$ such that 
	\[
	\vr \in L^\infty((0,T) \times \Td),\ \vr > 0,\ 
	\vu \in L^\infty((0,T) \times \Td; R^d). 
	\]
	
\end{Definition}

We are ready to state our main result. 

	\begin{mdframed}[style=MyFrame]
	
	\begin{Theorem} [\bf Density of wild data] \label{wT1}
		Suppose $p \in C^\infty(a,b)$, $p' > 0$ in $(a,b)$, for some $0 \leq a < b \leq \infty$. 
		
		Then for any
		\[
		\vr_0 \in W^{k,2}(\Td), \ a < \inf_{\Td} \vr_0 \leq \sup_{\Td} \vr_0 < b, 
		\vu_0 \in W^{k,2}(\Td; R^d),\ k > \frac{d}{2} + 1, 
		\]
		any $\ep > 0$, and any $1 \leq p < \infty$, there exist wild data $\vr_{0, \ep}$, $\vu_{0, \ep}$ such that 
		\[
		\| \vr_{0,\ep} - \vr_0 \|_{L^p(\Td)} < \ep,\ 
		\| \vu_{0, \ep} - \vu_0 \|_{L^p(\Td; R^d)} < \ep.
		\]		
		
	\end{Theorem}

\end{mdframed}

Recently, Chen, Vasseur, and You \cite{ChVaYu}, established density of wild data for the isentropic Euler system in the class of weak solutions satisfying the total energy inequality
\begin{equation} \label{ww1}
	\intTd{ \left[ \frac{1}{2} \vr |\vu|^2 + P(\vr) \right] (\tau, \cdot) } \leq \intTd{ \left[ \frac{1}{2} \vr_0 |\vu_0|^2 + P(\vr_0) \right] } \ \mbox{for any}\ \tau > 0.
	\end{equation}
These solutions are global in time, however, the associated total energy profile 
\[
\intTd{ \left[ \frac{1}{2} \vr |\vu|^2 + P(\vr) \right] (\tau, \cdot)  } 
\]
may not be non--increasing in time.

The proof of Theorem \ref{wT1} is based on a combination of the strong solution ansatz proposed by Chen, Vasseur, and You \cite{ChVaYu} with the abstract convex integration results 
concerning weak solutions with a given energy profile established in \cite{Fei2016}.

\section{Convex integration ansatz}
\label{C}

Similarly to Chen, Vasseur, and You \cite{ChVaYu}, our convex integration ansazt is based on strong solutions to the Euler system. 

\subsection{Local in time smooth solutions}

\begin{Proposition}[\bf Local existence for smooth data] \label{iP1}
	
	Suppose $p \in C^\infty(a,b)$, $p' > 0$ in $(a,b)$, for some $0 \leq a < b \leq \infty$. 
	
	Then for any initial data
	\begin{equation} \label{i9a}
	\vr_0 \in W^{k,2}(\Td), \ a < \inf_{\Td} \vr_0 \leq \sup_{\Td} \vr_0 < b, 
	\vu_0 \in W^{k,2}(\Td; R^d),\ k > \frac{d}{2} + 1
	\end{equation}
	there exists $T_{\rm max} > 0$ such that the compressible Euler system admits a classical solution 
	$(\vr, \vu)$ unique 
	in the class 
	\begin{equation} \label{i9}
	\vr \in C([0,T]; W^{k,2}(\Td)),\ a < \vr < b,\ 
	\vu \in C([0,T]; W^{k,2}(\Td; R^d))	
		\end{equation}
for any $0 < T < T_{\rm max}$.	
	\end{Proposition}

The proof of Proposition \ref{iP1} is nowadays standard and essentially attributed to Kato \cite{KatoT}, cf. also Benzoni-Gavage and Serre \cite[Chapter 13, Theorem 13.1] {BenSer}.

\subsection{Basic convex integration ansatz}
\label{C}

Consider the initial data $(\vr_0, \vu_0)$ in the regularity class \eqref{i9a} together with the associated 
smooth solution $(\tvr, \tvu)$ in $[0,T] \times \Td$, $T < T_{\rm max}$.
In addition, denote 
$\tvm = \tvr \tvu$. The Euler system written in the conservative variables $(\tvr, \tvm)$ reads 
\begin{align} 
	\partial_t \tvr + \Div \tvm &= 0, \label{C1} \\ 
	\partial_t \tvm + \Div \left( \frac{\tvm \otimes \tvm}{\tvr} + p(\tvr) \mathbb{I} \right) &= 0. 
	\label{C2}
	\end{align}

We look for solutions in the form 
\[
\vr = \tvr,\ \vm = \vr \vu = \tvm + \vv, 
\]
where 
\begin{align} 
	\Div \vv &= 0, \label{C3} \\ 
	\partial_t \vv + \Div \left( \frac{(\vv + \tvm) \otimes (\vv + \tvm) }{\tvr} - \frac{\tvm \otimes \tvm}{\tvr} \right) &=0, 
	\label{C4}, \\
	\vc{v}(0, \cdot) &= \vv_0. \label{C5}
	\end{align}

To apply the abstract results of \cite{Fei2016}, we rewrite problem \eqref{C3}--\eqref{C5} in the form
\begin{align} 
	\Div \vv &= 0, \label{C6} \\ 
	\partial_t \vv + \Div \left( \frac{(\vv + \tvm) \otimes (\vv + \tvm) }{\tvr} 
	- \frac{1}{d} \frac{|\vc{v} + \tvm |^2 }{\tvr} \mathbb{I}
	- \frac{\tvm \otimes \tvm}{\tvr}+ \frac{1}{d} \frac{|\tvm|^2}{\tvr} \mathbb{I} \right) &=0, 
	\label{C7} \\
	\vc{v}(0, \cdot) &= \vv_0. \label{C8}
\end{align}
\emph{together} with the prescribed ``kinetic energy''
\begin{equation} \label{C9}
	\frac{1}{2} \frac{|\vc{v} + \tvm|^2}{\tvr} = \frac{1}{2} \frac{|\tvm|^2}{\tvr} + \Lambda, 
\end{equation}
with a suitable spatially homogeneous function $\Lambda = \Lambda(t)$ determined below.

\section{Application of convex integration}
\label{a}

Setting 
\begin{equation} \label{a1}
\mathbb{H} = \frac{\tvm \otimes \tvm}{\tvr} - \frac{1}{d} \frac{|\tvm|^2}{\tvr} \mathbb{I} \in C^1([0,T] \times 
\Td; R^{d \times d}_{0, {\rm sym}}),\ e =  \frac{1}{2} \frac{|\tvm|^2}{\tvr} + \Lambda 
\in C([0,T] \times \Td),
\end{equation}
we may rewrite \eqref{C6}--\eqref{C9} in the form
\begin{align} 
	\Div \vv &= 0, \label{a2} \\ 
	\partial_t \vv + \Div \left( \frac{(\vv + \tvm) \otimes (\vv + \tvm) }{\tvr} 
	- \frac{1}{d} \frac{|\vc{v} + \tvm |^2 }{\tvr} \mathbb{I}
	- \mathbb{H} \right) &=0, 
	\label{a3} \\
	\frac{1}{2} \frac{|\vc{v} + \tvm|^2}{\tvr} &= e,  
	\label{a5} \\
	\vc{v}(0, \cdot) &= \vv_0, \label{a4}
\end{align}
for fixed $\mathbb{H}$ and $e$ given by \eqref{a1}.

\subsection{Subsolutions}

To apply the abstract results obtained in \cite{Fei2016}, 
we introduce the set of \emph{subsolutions} 
\begin{align} 
X_0 &= \left\{ \vv \in C_{\rm weak}([0,T]; L^2(\Td; R^d)) \cap L^\infty ((0,T) \times \Td; R^d) \ \Big| \ 
\right. \br	
&\quad \quad \vv(0,\cdot) = \vv_0,\ \vv(T; \cdot) = \vv_T, \vv \in C((0,T) \times \Td; R^d), \br 
&\quad \quad \Div \vv = 0,\ \partial_t \vv + \Div \mathbb{F} = 0 \ \mbox{in}\ \mathcal{D}'((0,T) \times \Td) \br
&\quad \quad \mbox{for some}\ \mathbb{F} \in L^\infty((0,T) \times \Td; R^{d \times d}_{0, {\rm sym}} ) 
\cap C((0,T) \times \Td; R^{d \times d}_{0, {\rm sym}} ),\br 
&\quad \quad \sup_{0 < \tau < t < T; x \in \Td } \frac{d}{2} \lambda_{\rm max} 
\left[ \frac{(\vv + \tvm) \otimes (\vv + \tvm) }{\tvr}     - \mathbb{F} - \mathbb{H} \right] - e < 0 , \br
&\quad \quad \mbox{for any}\ 0 < \tau < T \Big\}.
\label{a6}	
\end{align}
Here, the symbol $\lambda_{\rm max}[\mathbb{A}]$ denotes the maximal eigenvalue of a symmetric matrix $\mathbb{A}$.

\subsection{First existence result}

The following results is a special case of \cite[Theorem 13.2.1]{Fei2016}.

\begin{Proposition} \label{aP1}
	Suppose that set of subsolutions $X_0$ is non--empty and bounded in\\ $L^\infty((0,T) \times \Td; R^d)$, $d=2,3$. 
	
	Then problem \eqref{a2}--\eqref{a4} admits infinitely many weak solutions.
	
	\end{Proposition}

\bigskip

Fix $\vv_0 = \vv_T = 0$. Next, using the algebraic inequality
\begin{equation} \label{a7}
\frac{1}{2} \frac{ |\vv + \tvm|^2 }{\tvr} \leq 	
 \frac{d}{2} \lambda_{\rm max} 
\left[ \frac{(\vv + \tvm) \otimes (\vv + \tvm) }{\tvr}     - \mathbb{F} - \mathbb{H} \right]
\end{equation}
we can see that the set $X_0$ is bounded in $L^\infty((0,T) \times \Td; R^d)$ as long as $\Lambda \in C[0,T]$.
Finally, we observe that $\vv \equiv 0$ is a subsolution as soon as 
\begin{equation} \label{a8}
	\Lambda(t) > 0 \ \mbox{for any}\ t \in [0,T]. 
	\end{equation}
Indeed we may consider $\mathbb{F} \equiv 0$ and compute 
\[
\frac{d}{2} \lambda_{\rm max} 
\left[ \frac{(\vv + \tvm) \otimes (\vv + \tvm) }{\tvr}     - \mathbb{F} - \mathbb{H} \right] = 
\frac{d}{2} \lambda_{\rm max} 
\left[ \frac{\tvm \otimes \tvm }{\tvr}   - \mathbb{H} \right] = 
\frac{1}{2} \frac{|\tvm|^2}{\tvr}
\]
while 
\[
e = \frac{1}{2} \frac{|\tvm|^2}{\tvr} + \Lambda.
\] 

Thus a direct application of Proposition \ref{aP1} yields the following result. 

\begin{Theorem}[\bf Existence with a small initial energy jump] \label{aT1}
	
	Let $\Lambda \in C[0,T]$, $\inf_{t \in [0,T]} \Lambda(t) > 0$ be given. Let $\vv_0 = 0$. 
	
	Then problem \eqref{a1}--\eqref{a5} admits infinitely many weak solutions $\vv$ in $(0,T) \times \Td$.

	\end{Theorem}

As $\vv_0 = 0$ and $\Lambda > 0$, the solutions $\vr = \tvr$, $\vm = \tvm + \vv$ necessarily experience an initial energy jump therefore they are not physically admissible. This problem will be fixed in the next section.

\subsection{Second existence result}

The following results is a special case of \cite[Theorem 13.6.1]{Fei2016}.

\begin{Proposition} \label{aP2}
	
		Suppose that set of subsolutions $X_0$ is non--empty and bounded in\\ $L^\infty((0,T) \times \Td; R^d)$, $d=2,3$.
		
		Then there exists a set of time $\mathfrak{T} \subset (0,T)$ dense in $(0,T)$ with the following properties:
		
		For any $\tau \in \mathfrak{T}$ there exists $\vv^\tau \in \Ov{X}_0$ satisfying:
		
		\begin{itemize}
			\item 
			\begin{align} 
			\vv^\tau &\in C_{\rm weak}([0,T]; L^2(\Td; R^d)) \cap L^\infty ((0,T) \times \Td; R^d) \cap C((\tau, T) \times \Td; R^d),\ \br 
			\vv^\tau(0, \cdot) &= \vv_0,\ \vv^\tau(T, \cdot) = \vv_T; 	
				\label{a9}
				\end{align}
			
			\item 
			\begin{equation} \label{a10}
				\partial_t \vv^\tau + \Div \mathbb{F} = 0 \ \mbox{in}\ \mathcal{D}'((\tau, T) \times \Td) 
				\end{equation}
			for some $\mathbb{F} \in L^\infty \cap C((\tau; T) \times \Td; R^{d \times d}_{0, {\rm sym}})$;
			
			\item 
			\begin{equation} \label{a11}
			\sup_{\tau + s < t < T, x \in \Td} \frac{d}{2} \lambda_{\rm max} 
			\left[ \frac{(\vv^\tau + \tvm) \otimes (\vv^\tau + \tvm) }{\tvr}     - \mathbb{F} - \mathbb{H} \right] - e < 0	
				\end{equation}
			for any $0 < s < T - \tau$;
			
			\item
			
			\begin{equation} \label{a12}
			\frac{1}{2} \intTd{ \frac{|\vv^\tau + \tvm |^2}{\tvr}(\tau, \cdot) } = \intTd{ e(\tau, \cdot) }.
			\end{equation}
			
			\end{itemize} 
	
	\end{Proposition}

In accordance with \eqref{a9}--\eqref{a12}, the function $\vv^\tau$ can be used as a subsolution on the 
time interval $(\tau, T)$. Then Proposition \ref{aP1} yields the following result. 

\begin{Theorem} [\bf Existence without initial energy jump] \label{aT2}
	Let $\Lambda \in C[0,T]$, $\inf_{t \in [0,T]} \Lambda(t) > 0$ be given. 
		
	Then there exists a sequence $\tau_n \to 0$ and $\vv_{0,n}$, 
	\[
	\vv_{0,n} \to 0 \ \mbox{weakly-(*) in}\ L^\infty(\Td; R^d) 
	\]
	such that problem \eqref{a1}--\eqref{a5} admits infinitely many weak solutions in 
	$(\tau_n,T) \times \Td$ satisfying 
	\begin{equation} \label{a13}
	\vv(\tau_n, \cdot) = \vv_{0,n},\ \vv(T, \cdot) = 0,\ \frac{1}{2} \frac{|\vv + \tvm|^2}{\tvr} (\tau_n, \cdot) = e(\tau_n).	
		\end{equation}
	
	\end{Theorem}

Note carefully that continuity of the initial energy stated in \eqref{a13} follows from \eqref{a12} and weak continuity of $\vv$.

\section{Adjusting the energy profile}
\label{e}

To complete the proof of Theorem \ref{wT1}, it remains to adjust the energy profile $\Lambda$ so that: 

\begin{itemize} 
	\item
	\begin{equation} \label{e1}
	\limsup_{n \to 0} \| \vv_{0,n} \|_{L^2(\Td; R^d)} < \ep ;
		\end{equation}
	
	\item the energy inequality \eqref{i6} holds for 
	$\vu = \vv + \tvm$, $\vr = \tvr$, at least on a short time interval. 
	
	\end{itemize}
	
As for \eqref{e1}, it is enough to choose $\Lambda(0) > 0$ small enough. Indeed \eqref{a1}, \eqref{a13} 
yield 
\[
\frac{1}{2} \frac{|\vv + \tvm|^2}{\tvr} (\tau_n, \cdot) = \frac{1}{2} \frac{ |\tvm|^2 }{\tvr} (\tau_n) + 
\Lambda(\tau_n).
\]
Seeing that 
\[
\vv(\tau_n, \cdot) = \vv_{0,n} \to 0 \ \mbox{weakly in}\ L^2(\Td; R^d)
\]
we easily conclude.

Finally, the total energy of the system reads
\[
\frac{1}{2} \frac{ |\vv + \tvm |^2}{\tvr} + P(\tvr) = \frac{1}{2} \frac{|\tvm|^2}{\tvr} + \Lambda + P(\tvr)
\ \mbox{a.a. in}\ (0,T) \times \Td.
\]
In particular, the energy is continuously differentiable as soon as $\Lambda \in C^1[0,T]$. The desired energy 
inequality reads
\begin{equation} \label{e2}
\partial_t \left(\frac{1}{2} \frac{|\tvm|^2}{\tvr} + P(\tvr) \right) + \Lambda' + 
\Div \left[ \left(\frac{1}{2} \frac{|\tvm|^2}{\tvr} + P(\tvr) + \Lambda + p(\tvr) \right) \frac{ \tvm + \vv }{\tvr} \right] \leq 0.
\end{equation}
Seeing that the smooth solution $(\tvr, \tvm)$ satisfies the energy equality we may simplify \eqref{e2} to 
\begin{equation} \label{e3}
	\Lambda' + \Lambda \Div\tvu +
	\Div \left[ \left(\frac{1}{2} \frac{|\tvm|^2}{\tvr} + P(\tvr) + \Lambda + p(\tvr) \right) \frac{ \vv }{\tvr} \right] \leq 0.
\end{equation}
Moreover, as $\Div \vv = 0$,
\[ 
\Div \left[ \left(\frac{1}{2} \frac{|\tvm|^2}{\tvr} + P(\tvr) + \Lambda + p(\tvr) \right) \frac{ \vv }{\tvr} \right] 
= \Grad \left[ \frac{1}{\tvr} \left(\frac{1}{2} \frac{|\tvm|^2}{\tvr} + P(\tvr) + \Lambda + p(\tvr) \right) \right]
\cdot \vv, 
\]
and \eqref{e3}, reduces to 
\begin{equation} \label{e4}
	\Lambda' + \Lambda \Div\tvu +
\Grad \left[ \frac{1}{\tvr} \left(\frac{1}{2} \frac{|\tvm|^2}{\tvr} + P(\tvr) + \Lambda + p(\tvr) \right) \right]
\cdot \vv	\leq 0.
\end{equation}

As $\Lambda$ is decreasing in $t$, we get 
\[
\Lambda(t) \leq \Lambda(0). 
\]
Similarly, we control $\| \vv \|_{L^\infty((0,T) \times \Td; R^d)}$ by means of $\Lambda(0)$ and certain norms of the strong solution $\tvr$, $\tvm$.

Choosing 
\[
\Lambda(t) = \ep \exp \left( - \frac{t}{\ep^2} \right),
\]
with $\ep > 0$ small enough, we obtain the desired energy inequality at least on a short time interval $(0, T_{\rm w})$, $T_{\rm w} > 0$. We have proved Theorem \ref{wT1} for $p=2$. The same statement 
holds for a general $1 \leq p < \infty$ as all solutions in question  are uniformly bounded.

\section{Concluding remarks}
\label{cr}

A similar result can be shown for the more realistic complete slip boundary conditions 
\[
\vu \cdot \vc{n}|_{\partial \Omega} = 0
\]
imposed on a bounded domain $\Omega \subset R^d$. Note that the result is local in time and that the smooth solutions of the Euler system enjoy the finite speed of propagation property. Consequently, 
the problem of compatibility conditions may be solved by considering the initial data in the form 
\begin{align}
	\vr_0 \in W^{k,2}(\Td), \ a < \inf_{\Td} \vr_0 \leq \sup_{\Td} \vr_0 < b, 
	\vu_0 \in W^{k,2}(\Td; R^d),\ k > \frac{d}{2} + 1, \br
	\vu = 0, \vr = \Ov{\vr} - \mbox{a positive constant in a meighborhood of}\ \partial \Omega.
	\nonumber
\end{align}
The relevant local existence result for strong solutions was proved by Beir\~{a}o da Veiga \cite{BdaVei}.

\bibliography{citace}

\bibliographystyle{plain}

\end{document}